# Recent advances in invariance principles for stationary sequences*

**Florence Merlevède**

*Université Paris 6, LPMA and C.N.R.S UMR 7599, 175 rue du Chevaleret, 75013 Paris, FRANCE*
*e-mail:* merleve@ccr.jussieu.fr

**Magda Peligrad**†

*Department of Mathematical Sciences, University of Cincinnati, PO Box 210025, Cincinnati, Oh 45221-0025, USA*
*e-mail:* peligrm@math.uc.edu

**Sergey Utev**

*School of Mathematical Sciences, University of Nottingham, Nottingham, NG7 2RD, ENGLAND*
*e-mail:* sergey.utev@nottingham.ac.uk

**Abstract:** In this paper we survey some recent results on the central limit theorem and its weak invariance principle for stationary sequences. We also describe several maximal inequalities that are the main tool for obtaining the invariance principles, and also they have interest in themselves. The classes of dependent random variables considered will be martingale-like sequences, mixing sequences, linear processes, additive functionals of ergodic Markov chains.

## Contents

*This is an original survey paper

†Supported in part by a Charles Phelps Taft Memorial Fund grant and NSA grant H98230-05-1-0066.

## 1. Introduction

In the recent years, there has been a sustained effort towards a better understanding of the asymptotic behavior of stochastic processes.

For processes with short memory the theory of the weak invariance principle is very well fine tuned under various mixing conditions (see, e.g. the surveys by Peligrad [47], Philipp [58], Rio [61], Dehling and Philipp [14], and books by Bradley [5, 6]).

The only problem is that, in some situations, the mixing conditions are not verified. This is the reason that one of the new directions in modeling the dependence is to introduce new dependent structures defined by substantially reducing the classes of functions used in the definitions of mixing coefficients or by using innovative martingale-like conditions. In this way many new examples are included in the general structures and many general results can be established.

In this paper, we survey some recent results about the central limit theorem and its invariance principle under summability conditions imposed to the conditional expectation of consecutive sums of random variables with respect to the distant past.

For the sake of applications we shall survey also several results on the weak invariance principle under conditions imposed to the moments of the conditional expectation of the individual summand in a sequence. We shall see that these results extend the best results known for strongly mixing sequences and $\rho$–mixing sequences. At the same time, since the key technique for obtaining invariance principles is the use of maximal inequalities for partial sums, we shall also survey the recent advances on this subject.

In this paper we consider a stationary sequence $(X_k)_{k\in\mathbb{Z}}$ of centered random variables defined on a probability space $(\Omega, \mathcal{K}, \mathbb{P})$, with finite second moment and we shall survey some results on the central limit theorem and invariance principle namely

$$\frac{S_n}{b_n} \to \sqrt{\eta}N, \text{ as } n \to \infty, \tag{1}$$

where $S_n = \sum_{i=1}^n X_i$, $N$ has a standard normal distribution and $b_n \to \infty$, and

$$\frac{S_{[nt]}}{b_n} \to \sqrt{\eta} W(t) \text{ as } n \to \infty, \tag{2}$$

where $t \in [0,1]$, $W$ is the standard Brownian motion on $[0,1]$, and $\eta$ is a non-negative random variable. We shall discuss first several results for the situation when

$$\lim_{n\to\infty} \frac{b_n^2}{n} = c^2 \ .$$

We shall also point out some extensions of these theorems when the normalizer is more general, namely,

$$b_n^2 = \sigma_n^2 = \mathbb{E}(S_n^2) \qquad \text{or} \qquad b_n = \sqrt{\frac{2}{\pi}} \mathbb{E}|S_n| \ .$$

There is a vast literature on this subject. A certain restriction of the dependence structure is needed since a constant sequence obviously does not satisfy (1). Moreover also ergodicity or mixing in the ergodic sense are not sufficient. The classes of stochastic processes widely studied are martingales, uniformly mixing sequences, mixingales, associated sequences and so on.

## 2. Definitions and notations

Let $(\Omega, \mathcal{A}, \mathbb{P})$ be a probability space, and $\mathbb{T} : \Omega \mapsto \Omega$ be a bijective bi-measurable transformation preserving the probability. Let $\mathcal{F}_0$ be a $\sigma$-algebra of $\mathcal{A}$ satisfying $\mathcal{F}_0 \subseteq \mathbb{T}^{-1}(\mathcal{F}_0)$, and let $X_0$ be a random variable which is $\mathcal{F}_0$–measurable. We then define the nondecreasing filtration $(\mathcal{F}_i)_{i\in\mathbb{Z}}$ by $\mathcal{F}_i = \mathbb{T}^{-i}(\mathcal{F}_0)$ (referred to as the stationary filtration). We also define the stationary sequence $(X_i)_{i\in\mathbb{Z}}$ by $X_i = X_0 \circ \mathbb{T}^i$ (adapted to the stationary filtration $(\mathcal{F}_i)_{i\in\mathbb{Z}}$).

Alternatively, we may define a stationary filtration as in Maxwell and Woodroofe [35], that is we assume that $X_i = g(Y_j,\ j \le i)$ where $(Y_i)_{i\in\mathbb{Z}}$ is an underlying stationary sequence. Denote by $\mathcal{I}$ its invariant sigma field and by $(\mathcal{F}_i)_{i\in\mathbb{Z}}$ an increasing filtration of sigma fields $\mathcal{F}_i = \sigma(Y_j, j \le i)$. For the case when for every $i$, $\xi_i = Y_j$, and $g(Y_j, j \le i) = Y_i$, then $\mathcal{F}_i$ is simply the sigma algebra generated by $\xi_j$, $j \le i$.

We underline that we only consider a so-called adapted case, that is variables $X_i$ are $\mathcal{F}_i$–measurable.

Throughout the paper, we assume that $\mathbb{E}[X_0] = 0$ and $\mathbb{E}[X_0^2] < \infty$. Let

$$S_n = \sum_{k=1}^n X_k \ , \ \sigma_n^2 = \mathrm{Var}(S_n) \ , \ \|X\|_p = \mathbb{E}[|X|^p]^{1/p} \ .$$

Most of the results in this paper deal with stationary sequences, and sometimes the notion of invariant sigma field and ergodicity will be used. For all these definitions and canonical construction of stationary sequences we refer to the book of Bradley [5, Chapter 0, and Chapter 2].

**Definition 1** *We define the process* $\{W_n(t) : t \in [0,1]\}$ *by*

$$W_n(t) = \sum_{i=1}^{[nt]} X_i \ .$$

*For each* $\omega$, $W_n(\cdot)$ *is an element of the Skorohod space* $D([0,1])$ *of all functions on* $[0,1]$ *which have left-hand limits and are continuous from the right, equipped with the Skorohod topology (see, e.g. [1, Section 14]). Let* $W = \{W(t) \ ; \ t \in [0,1]\}$ *denote the standard Brownian motion on* $[0,1]$*. Throughout the paper by* $b_n^{-1} W_n \stackrel{\mathcal{D}}{\Longrightarrow} Z$ *we denote the convergence in distribution of the sequence of processes* $b_n^{-1} W_n = \{b_n^{-1} W_n(t) \ ; \ t \in [0,1]\}$ *in the Skorohod space* $D([0,1])$ *and by* $\stackrel{\mathcal{D}}{\to}$ *we denote the usual weak convergence on the real line.*

*Remark.* Alternatively, we could work with the process $Y_n(t) = S_{[nt]} + (nt - [nt])X_{[nt]+1}$, $t \in [0,1]$, which takes values in the space of continuous functions $C[0,1]$. We notice that in all our results, the limiting process $Z$ has continuous sample paths, since $Z$ has the form $Z = \sqrt{\eta} W$, where $W$ is the standard Brownian motion and variable $\eta$ is independent of $W$. Thus, all our results have an equivalent formulation as a weak convergence $b_n^{-1} Y_n \to Z$ in $C([0,1])$ (see, e.g. [1, Section 14]).

**Definition 2** *Following Definition 0.15 in Bradley [5], we will say that a sequence* $(h(n), n = 1, 2, 3, \dots)$ *of positive numbers is slowly varying in the strong sense (i.e. in the sense of Karamata) if there exists a continuous function* $f : (0, \infty) \to (0, \infty)$ *such that* $f(n) = h(n)$ *for all* $n \in N$*, and* $f(x)$ *is slowly varying as* $x \to \infty$*.*

We shall also introduce the following mixing coefficients: For any two $\sigma$-algebras $\mathcal{A}$ and $\mathcal{B}$ define the strong mixing coefficient $\alpha(\mathcal{A},\mathcal{B})$:

$$\alpha(\mathcal{A},\mathcal{B}) = \sup\{|\mathbb{P}(A \cap B) - \mathbb{P}(A)\mathbb{P}(B)|; \ A \in \mathcal{A} \ , \ B \in \mathcal{B}\} \ .$$

and the $\rho-$mixing coefficient, known also under the name of maximal coefficient of correlation $\rho(\mathcal{A},\mathcal{B})$:

$$\rho(\mathcal{A},\mathcal{B}) = \sup\{\mathrm{Cov}(X,Y)/\|X\|_2 \|Y\|_2 \ : \ X \in \mathbb{L}_2(\mathcal{A}), \ Y \in \mathbb{L}_2(\mathcal{B})\}.$$

(Here and in the sequel $\|.\|_2$ denotes the norm in $\mathbb{L}_2$.)

For the stationary sequence of random variables $(X_k)_{k \in \mathbb{Z}}$, we also define $\mathcal{F}_m^n$ the $\sigma$–field generated by $X_i$ with indices $m \le i \le n$, $\mathcal{F}^n$ denotes the $\sigma$–field generated by $X_i$ with indices $i \ge n$, and $\mathcal{F}_m$ denotes the $\sigma$–field generated by $X_i$ with indices $i \le m$. Notice that $(\mathcal{F}_k)_{k \in \mathbb{Z}}$ defined in this way is a minimal filtration. The sequences of coefficients $\alpha(n)$ and $\rho(n)$ are then defined by

$$\alpha(n) = \alpha(\mathcal{F}_0, \mathcal{F}^n), \text{ and } \ \rho(n) = \rho(\mathcal{F}_0, \mathcal{F}^n) \ .$$

Equivalently, (see [5, ch. 4])

$$\rho(n) = \sup\{\|\mathbb{E}(Y|\mathcal{F}_0)\|_2 / \|Y\|_2 \ : \ Y \in \mathbb{L}_2(\mathcal{F}^n), \ \mathbb{E}(Y) = 0\} \ . \tag{3}$$

Finally we say that the stationary sequence is strongly mixing if $\alpha(n) \to 0$ as $n \to \infty$, and $\rho-$mixing if $\rho(n) \to 0$ as $n \to \infty$.

In some situations we are going to use weaker forms of strong and $\rho$–mixing coefficients, when $\mathcal{F}^n$, is replaced by the sigma algebra generated by only one variable, $X_n$, denoted by $\mathcal{F}_n^n$. We shall use the notations $\tilde{\alpha}(n) = \alpha(\mathcal{F}_0,\mathcal{F}_n^n)$ and $\tilde{\rho}(n) = \rho(\mathcal{F}_0,\mathcal{F}_n^n)$.

## 3. Maximal inequalities

In this section, unless otherwise specified, we do not assume stationarity. It is convenient to use the same notation, $(X_i)_{i\in\mathbb{Z}}$, for a sequence of square-integrable random variables, not necessarily stationary, adapted to the nondecreasing filtration $(\mathcal{F}_i)_{i\in\mathbb{Z}}$, also non-stationary, in general. Let

$$\begin{aligned} S_n &= \sum_{k=1}^{n} X_k, \\ M_n &= \max_{1\le i\le n} |S_i|, \\ M_n^+ &= \max_{1\le j\le n} S_j^+ = \max(0, S_1, ..., S_n) \text{ and} \\ M_n^- &= \max_{1\le j\le n} (-S_j^-) = \max(0, -S_1, ..., -S_n) \,. \end{aligned}$$

We shall start this section with a result which is a common step for obtaining several maximal inequalities. Next result was formulated in Peligrad and Utev [55, relation (10)]. Its proof is based on the interesting inequality suggested in Dedecker and Rio [12, relation (3.4)], which was obtained by using Garsia's [19] telescoping sums approach to the maximal inequality.

**Proposition 3** *Let $(X_i)_{i\in\mathbb{Z}}$ be a sequence of random variables. Then*

$$(M_n)^2 \le 4(S_n)^2 - 4\sum_{k=1}^{n-1} D_k(S_n - S_k) \,, \tag{4}$$

*where* $D_k = (M_k^+ - M_{k-1}^+) - (M_k^- - M_{k-1}^-)$.

By taking the expected values in (4), we then derive the following result:

**Proposition 4** *Let $(X_i)_{i\in\mathbb{Z}}$ be a sequence of square-integrable random variables, adapted to a nondecreasing filtration $(\mathcal{F}_i)_{i\in\mathbb{Z}}$. Then*

$$\mathbb{E}(M_n)^2 \le 4\mathbb{E}(S_n)^2 + 4\sum_{k=1}^{n-1} \|X_k \mathbb{E}(S_n - S_k|\mathcal{F}_k)\|_1 \,. \tag{5}$$

Dedecker and Rio [12] further expanded $\mathbb{E}(S_n^2)$ in (5) and obtained the following extension of Doob's inequality for martingales.

**Proposition 5** *Let $(X_i)_{i\in\mathbb{Z}}$ be a sequence of square-integrable random variables, adapted to a nondecreasing filtration $(\mathcal{F}_i)_{i\in\mathbb{Z}}$. Let $\lambda$ be any nonnegative real number and $\Gamma_k = (M_k^+ > \lambda)$.*
*(a) We have*

$$\begin{aligned}\mathbb{E}((M_n^+ - \lambda)_+^2) \;\leq\; & 4\sum_{k=1}^{n} \mathbb{E}(X_k^2 \; \mathbb{1}_{\Gamma_k}) \\ & +8\sum_{k=1}^{n-1} ||X_k \; \mathbb{1}_{\Gamma_k} \mathbb{E}(S_n - S_{k-1}|\mathcal{F}_k)||_1 \; .\end{aligned}$$

*(b) If furthermore the two-dimensional array $(X_k\mathbb{E}(S_n - S_{k-1}|\mathcal{F}_k))_{1\leq k\leq n}$ is uniformly integrable then the sequence $(n^{-1}M_n^2)_{n>0}$ is uniformly integrable.*

We want also to point out a somehow less general version of it but with better constants, which is easy to derive from Proposition 4.

**Corollary 6** *Let $(X_i)_{i\in\mathbb{Z}}$ be a sequence of square-integrable random variables, adapted to a nondecreasing filtration $(\mathcal{F}_i)_{i\in\mathbb{Z}}$. Then*

$$\mathbb{E}(M_n)^2 \leq 4\sum_{k=1}^{n} \mathbb{E}(X_k^2) + 12\sum_{k=1}^{n-1} \|X_k\mathbb{E}(S_n - S_k|\mathcal{F}_k)\|_1 \; . \tag{6}$$

Concerning the higher moments, we give the following result which is an extension to non-stationary sequences of Theorem 2.5 of Rio [61] and which is due to Dedecker and Doukhan [9].

**Proposition 7** *Let $(X_i)_{i\in\,\mathbf{N}}$ be a sequence of square integrable random variables, and $\mathcal{F}_i = \sigma(X_j, 0 \leq j \leq i)$. Define $S_n = X_1 + \cdots + X_n$ and*

$$b_{i,n} = \max_{i\leq l\leq n} \|X_i \sum_{k=i}^{l} \; \mathbb{E}(X_k|\mathcal{F}_i)\|_{p/2} \; .$$

*For any $p \geq 2$, the following inequality holds*

$$\mathbb{E}|S_n|^p \leq (2p\sum_{i=1}^{n} b_{i,n})^{p/2} \; . \tag{7}$$

*In addition for $p \geq 2$,*

$$\mathbb{E}(\max_{1\leq j\leq n} |S_j|^p) \leq C_p \, (\sum_{i=1}^{n} b_{i,n})^{p/2} \; . \tag{8}$$

For $p = 2$, (8) follows from (6) with $C_p = 16$. For $p > 2$, (8) is a direct consequence of (7) combined with Theorem 1 in Móricz [44]. In this case $C_p = (1 - 2^{(1-p)/(2p)})^{-2p}(2p)^{p/2}$.

We now state a maximal inequality which was derived from an adaptation of a result by McLeish [38] in combination with Proposition 5(a) (see Proposition 6 in Dedecker and Merlevède [10]). To describe this result, we consider the projection operator $P_i$:

$$\text{for any } f \in \mathbb{L}_2,\ P_i(f) = \mathbb{E}(f|\mathcal{F}_i) - \mathbb{E}(f|\mathcal{F}_{i-1}) \,. \tag{9}$$

**Proposition 8** *Let $(X_i)_{i\in\mathbb{Z}}$ be a sequence of square–integrable random variables, adapted to a nondecreasing filtration $(\mathcal{F}_i)_{i\in\mathbb{Z}}$. Define the $\sigma$–algebra $\mathcal{F}_{-\infty} = \bigcap_{i\in\mathbb{Z}} \mathcal{F}_i$. For any positive integer $i$, let $(Y_{i,j})_{j\geq 1}$ be the martingale defined by*

$$Y_{i,j} = \sum_{k=1}^{j} P_{k-i}(X_k) \quad \text{and let} \quad Y_{i,n}^{+} = \max\{0, Y_{i,1}, \ldots, Y_{i,n}\} \,.$$

*Assume that for any integer $k$, $\mathbb{E}(X_k|\mathcal{F}_{-\infty}) = 0$. For $\lambda \geq 0$ and nonnegative sequences $(a_i)_{i\geq 0}$ and $(b_i)_{i\geq 0}$ with $K = \sum a_i^{-1} < \infty$ and $\sum b_i = 1$, we have*

$$\mathbb{E}\left((M_n^{+} - \lambda)_{+}^{2}\right) \leq 4K \sum_{i=0}^{\infty} a_i \Big( \sum_{k=1}^{n} \mathbb{E}(P_{k-i}^{2}(X_k)\ \mathbb{1}_{(Y_{i,k}^{+} > b_i\lambda)}) \Big) \,. \tag{10}$$

We now move to present a new maximal inequality for stationary sequences due to Peligrad and Utev [55]. With this aim, we introduce the notation:

$$\Delta_r = \sum_{j=0}^{r-1} \left\| \frac{\mathbb{E}(S_{2^j}|\mathcal{F}_0)}{2^{j/2}} \right\|_2 \,. \tag{11}$$

**Proposition 9** *Let $(X_i)_{i\in\mathbb{Z}}$ be a stationary sequence of random variables with $\mathbb{E}(X_0) = 0$ and $\mathbb{E}(X_0^2) < \infty$. Let $n$, $r$ be integers such that $2^{r-1} < n \leq 2^r$. Then we have*

$$\mathbb{E}\left[\max_{1\leq i\leq n} S_i^2\right] \leq n\left(2\|X_1\|_2 + (1+\sqrt{2})\Delta_r\right)^2 \,. \tag{12}$$

The proof of this maximal inequality is based on several ideas including:

($i$) the diadic induction found to be useful in the analysis of $\rho$–mixing sequences (see [46], [8], [54], [63] for other variations on the theme),

($ii$) the modification of the Garsia [19] telescoping sums approach to maximal inequalities as in Peligrad [52] and Dedecker and Rio [12] and

($iii$) the subadditivity properties of the conditional expectation of sums.

To comment on the quantity defined in relation (11), we notice that, by the proof of Lemma 3.3 in Peligrad and Utev [55], we can find two positive constants $C_1$ and $C_2$ such that for $2^{r-1} < n \leq 2^r$

$$C_1 \sum_{j=1}^{n} \left\| \frac{\mathbb{E}(S_j|\mathcal{F}_0)}{j^{3/2}} \right\|_2 \leq \Delta_r \leq C_2 \sum_{j=1}^{n} \left\| \frac{\mathbb{E}(S_j|\mathcal{F}_0)}{j^{3/2}} \right\|_2 \,.$$

In a recent paper Peligrad, Utev and Wu [57] extended the inequality (12) to moments higher than 2 and obtained also an exponential bound. Basically, an estimate for $\mathbb{E}\left[\max_{1\le i\le n}|S_i|^p\right]$ with $p\ge 2$ was obtained in terms of $\|X_1\|_p$ and $\sum_{j=1}^n j^{-3/2}\left\|\mathbb{E}(S_j|\mathcal{F}_0)\right\|_p$.

Of course, there are many other useful maximal inequalities for dependent structures. We mention a recent survey by Bradley [4] that contains useful connections between mixing properties, inequalities and normal approximation type results.

As an example of other possible approaches, we state the following powerful Rosenthal type moment maximal inequality for sums of random variables in terms of the interlaced mixing coefficients.

For the stationary sequence $\{X_k\}_{k\in\mathbb{Z}}$ denote by $\mathcal{F}_T=\sigma(X_i; i\in T)$ where $T$ is a family of integers. Define

$$\rho_n^* = \sup \rho(\mathcal{F}_T, \mathcal{F}_S)$$

where this supremum is taken over all pairs of nonempty finite sets $S, T$ of $Z$ such that $\mathrm{dist}(S,T)\ge n$. The following inequality follows. See Peligrad [50, 51], Peligrad and Gut [53], Utev and Peligrad [66] for details.

**Theorem 10** *For a positive integer $N\ge 1$ and positive real numbers $q\ge 2$ and $0\le r<1$, there is a positive constant $D=D(q,N,r)$ such that if $\{X_i\}_{i\ge 1}$ is a sequence of random variables with $\rho_N^*\le r$, and such that $\mathbb{E}X_i=0$ and $\mathbb{E}|X_i|^q<\infty$ for every $i\ge 1$, then for all $n\ge 1$,*

$$\mathbb{E}\Big(\max_{1\le i\le n}\Big|\sum_{j=1}^{i}X_j\Big|^q\Big)\le D\left(\sum_{i=1}^n \mathbb{E}|X_i|^q+\left(\sum_{i=1}^n \mathbb{E}X_i^2\right)^{\frac{q}{2}}\right) .$$

## 4. Invariance principles

In this section we survey several recent results on the central limit and its invariance principles for martingale-like sequences. The first part deals with the normalization $\sqrt{n}$.

### 4.1. The normalization $\sqrt{n}$

One of the possible approaches to study the asymptotic behavior of the normalized partial sum process is to approximate $S_n$ by a related martingale with stationary differences. Then, under some additional conditions, the central limit theorem can be deduced from the martingale case. This approach was first explored by Gordin [20], who obtained a sufficient condition for the asymptotic normality of the normalized partial sums. To explain the method, we shall introduce the perturbation approach, motivated by Liverani [34].

For a suitable sequence $(a_n)_{n\geq 1}$, define the sequence of random variables:

$$\theta_k = \sum_{i=0}^{\infty} a_i \mathbb{E}(X_{k+i}|\mathcal{F}_k) \ .$$

We shall assume that $(a_n)_{n\geq 1}$ and $(X_i)_{i\in\mathbb{Z}}$ satisfy conditions that assure that the above series is convergent in $\mathbb{L}_1$ and $\theta_k$ denotes the limit in $\mathbb{L}_1$ of this series. Natural conditions to assume are: $\mathbb{E}|X_0| < \infty$ , $\sum_{i=0}^{\infty} |a_i| < \infty$ and $a_0 = 1$, but other situations are also possible. We also notice that

$$\theta_k = X_k + \sum_{i=1}^{\infty} a_i \mathbb{E}(X_{k+i}|\mathcal{F}_k)$$

and

$$\mathbb{E}(\theta_{k+1}|\mathcal{F}_k) - \theta_k = -X_k + \sum_{i=1}^{\infty} (a_{i-1} - a_i) \mathbb{E}(X_{k+i}|\mathcal{F}_k) \ . \tag{13}$$

Finally, denote by

$$D_k = \theta_{k+1} - \mathbb{E}(\theta_{k+1}|\mathcal{F}_k) \ .$$

Then, $(D_k)_{k\in\mathbb{Z}}$ is a martingale difference sequence which is stationary (and ergodic) if the sequence $(X_i)_{i\in\mathbb{Z}}$ is stationary with the stationary filtration $(\mathcal{F}_i)_{i\in\mathbb{Z}}$ (and ergodic). By summing over $k$ in (13) we obtain the following form of the so-called coboundary decomposition

$$S_n = M_n + \sum_{k=1}^{n} \sum_{i=1}^{\infty} (a_{i-1} - a_i) \mathbb{E}(X_{k+i}|\mathcal{F}_k) + \theta_1 - \theta_{n+1} \tag{14}$$

where $M_n = \sum_{k=1}^{n} D_k$ is a martingale. Now, if we are interested for instance in proving the CLT for $S_n$ we have just to apply the well-known CLT theorem to the ergodic martingale $M_n$ and therefore we have just to impose conditions that give

$$\frac{1}{\sqrt{n}} \mathbb{E} \left| \sum_{k=1}^{n} \sum_{i=1}^{\infty} (a_{i-1} - a_i) \mathbb{E}(X_{k+i}|\mathcal{F}_k) \right| \to 0 \text{ as } n \to \infty \ . \tag{15}$$

In the martingale approach to stationary sequences (Gordin [20]), the selection of $(a_n)_{n\geq 1}$ is $a_n = 1$ for all $n$ and also the key condition is roughly,

$$\mathbb{E}(S_n|\mathcal{F}_0) \text{ is convergent in } \mathbb{L}_2 \ . \tag{16}$$

This condition was further relaxed by Gordin [21], Hall and Heyde [22], Peligrad [45], Volný [67] and others to

$$\mathbb{E}(S_n|\mathcal{F}_0) \text{ is convergent in } \mathbb{L}_1 \ . \tag{17}$$

The conditions (16) and (17) are motivated by the so called coboundary problem. More details can be found in the book by Hall and Heyde [22] and the correction by Esseen and Janson [18]. The following result was announced by Gordin and Lifsic at the third Vilnius conference (1981): the central limit theorem holds if $P_0(S_n)$ is convergent in $L_2$, where $P_0$ is the projection operator (defined by (9)) and $\mathbb{E}(S_n|\mathcal{F}_0)/\sqrt{n}$ converges to 0 in $L_2$ (or eventually in $L_1$), ( see, Borodin and Ibragimov [2, Theorem 8.1] for more details). For some related papers on Markov chains we mention papers by Kipnis and Varadhan [33], Woodroofe [68], Derriennic and Lin [15].

Maxwell and Woodroofe [35] and Wu and Woodroofe [70] considered the selection of the numerical sequence $a_k = a_k(m) = (1+1/m)^{-k-1}$ , where $m$ is a parameter that is allowed to depend on $n$ when a partial sum of size $n$ is considered. This selection together with the condition

$$\sum_{n=1}^{\infty} \frac{\|\mathbb{E}(S_n|\mathcal{F}_0)\|_2}{n^{3/2}} < \infty \tag{18}$$

assure the validity of the martingale decomposition above. They proved the central limit theorem under the assumption (18).

Now, a natural question is to derive sufficient conditions for the validity of the invariance principle. Considering again the perturbation approach, by replacing $n$ with $[nt]$ in (14), we only have to apply the classical invariance principle to the ergodic martingale $\frac{M_{[nt]}}{\sqrt{n}}$ and prove that, for every $\varepsilon > 0$,

$$\lim_{n\to\infty} \mathbb{P}\left( \max_{1\le j\le n} \left| \sum_{k=1}^{j} \sum_{i=1}^{\infty} (a_{i-1} - a_i)\mathbb{E}(X_{k+i}|\mathcal{F}_k) \right| \ge \varepsilon\sqrt{n} \right) = 0\,, \tag{19}$$

and

$$\lim_{n\to\infty} \mathbb{P}\left( \max_{1\le j\le n} |\theta_j| \ge \varepsilon\sqrt{n} \right) = 0\,. \tag{20}$$

Then, selecting $(a_n)_{n\ge1}$ such that $a_n = 1$ for all $n$, it obviously remains to prove (20) which follows from standard computations under (16) (see also Heyde [25, 26] and Volný [67]). However, criterion (16) may be suboptimal when applied to Markov chains or to strongly mixing sequences. Now, it seems quite natural to answer if the invariance principle holds under the weaker condition (17). Nevertheless, it appears from Remark 3 in Volný [67], that (17) is not a sufficient condition to get the weak invariance principle. However, it is possible to weaken (16). Indeed, recently Peligrad and Utev [55] proved the following invariance principle under the condition (18):

**Theorem 11** *Let $(X_i)_{i\in\mathbb{Z}}$ be a stationary sequence with $\mathbb{E}(X_0) = 0$ and $\mathbb{E}(X_0^2) < \infty$. Assume that (18) holds. Then,*

$$\{\max_{1\le k\le n} S_k^2/n \;:\; n\ge 1\} \quad \textit{is uniformly integrable and} \quad n^{-1/2}W_n \stackrel{\mathcal{D}}{\Longrightarrow} \sqrt{\eta}W\,,$$

*where $\eta$ is a non-negative random variable with finite mean $\mathbb{E}[\eta] = \sigma^2$ and independent of $\{W(t); t\ge 0\}$. Moreover, $\eta$ is determined by the limit*

$$\lim_{n\to\infty} \frac{\mathbb{E}(S_n^2|\mathcal{I})}{n} \;=\; \eta \;\; \textit{in } \mathbb{L}_1\,, \tag{21}$$

*where $\mathcal{I}$ is the invariant sigma field. In particular,* $\lim_{n\to\infty}\mathbb{E}(S_n^2)/n=\sigma^2$.

For the sake of applications, our next corollary formulates a sufficient condition for the invariance principles in terms of the conditional expectation of an individual summand $X_n$ with respect to $\mathcal{F}_0$. The following result strengthens Corollary 5 in Maxwell and Woodroofe [35].

**Corollary 12** *Assume that*

$$\sum_{n=1}^{\infty}\frac{1}{\sqrt{n}}\|\mathbb{E}(X_n|\mathcal{F}_0)\|_2\ <\ \infty\ . \tag{22}$$

*Then (18) is satisfied, and the conclusion of Theorem 11 holds.*

Let us notice that the condition (22) entails (18) by using the relation

$$\|\mathbb{E}(S_n|\mathcal{F}_0)\|_2\le\sum_{j=1}^{n}\|\mathbb{E}(X_j|\mathcal{F}_0)\|_2\ .$$

By Hölder inequality, (22) is satisfied as soon as there exists a sequence $(L_n)_{n\ge1}$ of positive numbers such that:

$$\sum_{n=1}^{\infty}\frac{1}{nL_n}\ <\ \infty\quad\text{and}\quad\sum_{n=1}^{\infty}L_n\|\mathbb{E}(X_n|\mathcal{F}_0)\|_2^2<\infty\ . \tag{23}$$

Note also that (22) entails (23) by simply taking $L_n=\big(\sqrt{n}\|\mathbb{E}(X_n|\mathcal{F}_0)\|_2\big)^{-1}$. Consequently (22) and (23) are equivalent.

Then, by using (23), the conclusion of Theorem 11 holds under

$$\sum_{n=1}^{\infty}\log^v(n)\|\mathbb{E}(X_n|\mathcal{F}_0)\|_2^2\ <\ \infty,\ \text{for some } v>1\ . \tag{24}$$

**Idea of the proof of Theorem 11.** In order to get the invariance principle (i.e. the convergence of the stochastic processes $n^{-1/2}W_n(t)$ to the Brownian motion) and at the same time to avoid the requirement of ergodicity, Peligrad and Utev [55] followed a different approach that is described in what follows. The initial step of the proof was a decomposition in small blocks of size $m$,

$$X_i^{(m)}=\frac{1}{\sqrt{m}}\sum_{j=(i-1)m+1}^{im}X_j\ .$$

Then a sequence of martingale differences was constructed,

$$D_i^{(m)}=X_i^{(m)}-\mathbb{E}(X_i^{(m)}|\mathcal{F}_{i-1}^{(m)})$$

(where $\mathcal{F}_\ell^{(m)}$ is a $\sigma$-field generated by variables $X_j^{(m)}$ with indices $j\le\ell$) and setting $k=[n/m]$ (where, as before, $[x]$ denotes the integer part of $x$), the following stationary martingale is introduced:

$$M_k^{(m)} = \sum_{i=1}^{k} D_i^{(m)} \ .$$

By using the classical invariance principle for martingales, it follows that

$$\frac{1}{\sqrt{k}} M_{[kt]}^{(m)} \overset{\mathcal{D}}{\Longrightarrow} \sqrt{\eta^{(m)}} W \ .$$

where $\eta^{(m)}$ is the limit (both : in $\mathbb{L}_1$ and almost surely) of

$$\frac{1}{k} \sum_{i=1}^{k} \Big( X_i^{(m)} - \mathbb{E}(X_i^{(m)} | \mathcal{F}_{i-1}^{(m)}) \Big)^2 \qquad \text{as} \quad k \to \infty \ .$$

In order to prove the invariance principle for $\frac{1}{\sqrt{n}} S_{[nt]}$, together with the uniform integrability of the sequence

$$\max_{1 \le k \le n} S_k^2 / n \ ,$$

it was established that

$$\| \sqrt{\eta^{(m)}} - \sqrt{\eta} \|_2 \to 0 \text{ as } m \to \infty \tag{25}$$

and

$$\lim_{m \to \infty} \lim_{n \to \infty} \left\| \sup_{0 \le t \le 1} \left| \frac{1}{\sqrt{n}} S_{[nt]} - \frac{1}{\sqrt{k}} M_{[kt]}^{(m)} \right| \right\|_2 = 0 \ . \tag{26}$$

The main part of the proof was to show (26). After some computations, this was reduced to establishing that

$$\lim_{m \to \infty} \lim_{n \to \infty} \frac{1}{\sqrt{k}} \left\| \max_{1 \le j \le k} \left| \sum_{i=1}^{j} \mathbb{E}(X_i^{(m)} | \mathcal{F}_{i-1}^{(m)}) \right| \right\|_2 = 0 \,, \tag{27}$$

that was derived by applying the maximal inequality from Proposition 9.

As a counterpart to Theorem 11, we now state the invariance principle derived in Dedecker and Rio [12].

**Theorem 13** *Let $(X_i)_{i \in \mathbb{Z}}$ be a stationary sequence with $\mathbb{E}(X_0) = 0$ and $\mathbb{E}(X_0^2) < \infty$. Assume that*

$$\mathbb{E}(X_0 S_n | \mathcal{F}_0) \ \textit{converges in} \ \mathbb{L}_1 \ . \tag{28}$$

*Then $n^{-1/2} W_n \overset{\mathcal{D}}{\Longrightarrow} \sqrt{\eta} W$ , where $W$ is a standard Brownian motion independent of $\mathcal{I}$ and $\eta$ is the $\mathcal{I}$-measurable non-negative variable which is determined by the limit $\lim_{n \to \infty} (\mathbb{E}(X_0^2 | \mathcal{I}) + 2\mathbb{E}(X_0 S_n | \mathcal{I}))_{n>0} = \eta$ in $\mathbb{L}_1$.*

**Remark 1.** If the sequence is ergodic, (that is the sigma field $\mathcal{I}$ is trivial) then

$$\eta = \sigma^2 = \mathbb{E}(X_0^2) + 2\sum_{i>0} \mathbb{E}(X_0 X_i)$$

and the usual invariance principle holds.

The following corollary was frequently used in various examples in Dedecker and Rio [12].

**Corollary 14** *Assume that*

$$\sum_{i=1}^{\infty} \mathbb{E}|X_0 \mathbb{E}(X_i|\mathcal{F}_0)| < \infty \ . \tag{29}$$

*Then the conclusion of Theorem 13 holds.*

**On the proof of Theorem 13**. The proof of this theorem involves a delicate application of Lindeberg approach together with the usual blocking technique, and the use of Proposition 5. It would be interesting to see whether a similar approach to the one used in the proof of Theorem 11 could work in this situation. The proof could follow the same line except for, in order to prove (27), Corollary 6 should be used instead of Proposition 9. However this approach apparently requires very delicate computations. With this possible approach in mind we raise here the following questions:

For example, assume (29). Does it follow that

$$\sum_{i=1}^{\infty} \mathbb{E}\left|X_0^{(m)} \mathbb{E}(X_i^{(m)}|\mathcal{F}_0)\right| \ \ \to \ \ 0, \quad \text{as } m \to \infty \,?$$

Also assume that (28) holds. Does it follow that

$$\lim_{m\to\infty} \limsup_{n\to\infty} \mathbb{E}\left|X_0^{(m)} \sum_{i=1}^{n} \mathbb{E}(X_i^{(m)}|\mathcal{F}_0)\right| = 0\,?$$

### *4.2. Examples. Part 1*

In this section, we provide some examples where the conditions are written in more familiar terms specialized to examples considered. Whenever possible, we will compare performance of Theorems 11, and 13 , and their Corollaries 12 and 14 for various examples.

### 4.2.1. Mixing sequences

In this group of examples, we apply the general results to stationary mixing sequences. Starting with the seminal works by Rosenblatt [62] and Ibragimov [27], the central limit theorem, maximal inequalities and the invariance principle were thoroughly analyzed by many mathematicians. In some mixing cases such as $\alpha$-mixing (strong mixing) and $\rho$ -mixing, optimal sufficient conditions for the validity of the CLT and its invariance principles were found, supported by various counterexamples. Yet, the theory is far from being completely understood. Various examples of processes satisfying mixing conditions can be found in the recent books by Doukhan [16] and by Bradley [5, 6].

This class of examples will show how the results surveyed in the previous section can be compared to known results in the literature using standard mixing type methods.

First, we point that all examples considered in this part satisfy

$$\mathbb{E}|\mathbb{E}(X_n|\mathcal{F}_0)|^2 = \|\mathbb{E}(X_n|\mathcal{F}_0)\|_2^2 \to 0 \text{ as } n \to \infty \ . \tag{30}$$

We notice that condition (29) is practically more restrictive and demands a certain rate of convergence in (30).

*a) $\rho$-mixing sequences.*

An interesting point of Theorem 11 is its parallel with the maximal coefficient of correlation $\rho$ associated to the stationary sequence $(X_k)_{k\in\mathbb{Z}}$. It is well known that the central limit theorem and its invariance principle hold under the assumption

$$\sum_{k=1}^{\infty} \rho(2^k) < \infty \ . \tag{31}$$

where $\rho(n) = \rho(\mathcal{F}_0,\mathcal{F}^n)$. Let us recall that the central limit theorem is due to Ibragimov [30], but the invariance principle remained an open question for a long time and it took an intensive effort from the part of many mathematicians, including Peligrad [46], Shao [63] and Utev [64, 65] to solve this problem. In addition, we would like to mention that the class of processes satisfying the condition (18) is larger than the class of function satisfying (31).

We shall see that for $\rho-$mixing sequences Theorem 11 gives optimal results, while for this case, a direct application of Corollary 14 will require a much more restrictive polynomial rate, namely $\sum_{k=1}^{\infty} \rho(k) < \infty$ that has been considered as a best possible for about ten years.

To prove this and similar results, the following identity taken from Dedecker and Rio [12, (6.1)], appears to be useful

$$\mathbb{E}|X\mathbb{E}(Y|\mathcal{F})| = \mathrm{Cov}(|X|(\ \mathbb{1}_{\mathbb{E}(Y|\mathcal{F})>0} - \ \mathbb{1}_{\mathbb{E}(Y|\mathcal{F})\le 0}), Y) \ , \tag{32}$$

where $X$ is $\mathcal{F}-$measurable and $Y$ is centered.

In particular, by Corollary 14 and (32) applied with $X = X_0$, $Y = X_k$ and $\mathcal{F} = \mathcal{F}_0$, it follows that a weaker condition, namely $\sum_{k=1}^{\infty} \tilde{\rho}(k) < \infty$ , is sufficient for the central limit theorem and invariance principle. Moreover, Corollary 12 needs even a weaker assumption

$$\sum_{k=1}^{\infty} \frac{\tilde{\rho}(k)}{\sqrt{k}} < \infty \,.$$

Let us now give a short argument to prove that (31) implies (18). We mention first the well-known fact that (31) implies,

$$\mathbb{E}[S_n^2] \le c^2 n \qquad (c>0)\,. \tag{33}$$

By stationarity we have:

$$\|\mathbb{E}(S_{2n}|\mathcal{F}_0)\|_2 \le \|\mathbb{E}(S_n|\mathcal{F}_0)\|_2 + \|\mathbb{E}(S_n|\mathcal{F}_{-n})\|_2 \,.$$

Clearly,

$$\mathbb{E}[\mathbb{E}(S_n|\mathcal{F}_{-n})^2] = \mathrm{Cov}(S_n, \mathbb{E}(S_n|\mathcal{F}_{-n})) \le \rho(n)\|S_n\|_2 * \|\mathbb{E}(S_n|\mathcal{F}_{-n})\|_2 \tag{34}$$

so that $\|\mathbb{E}(S_n|\mathcal{F}_{-n})\|_2 \le c\rho(n)\sqrt{n}$. Thus by recurrence,

$$\|\mathbb{E}(S_{2^{r+1}}|\mathcal{F}_0)\|_2 \le \|\mathbb{E}(S_{2^r}|\mathcal{F}_0)\|_2 + c2^{r/2}\rho(2^r) \le \ldots \le c\sum_{j=0}^{r} 2^{j/2}\rho(2^j)$$

and

$$\sum_{r=0}^{\infty} \frac{\|\mathbb{E}(S_{2^r}|\mathcal{F}_0)\|_2}{2^{r/2}} \;\le\; c\sum_{j=0}^{\infty} 2^{j/2}\rho(2^j)\sum_{r=j}^{\infty} 2^{-r/2} \le 4c\sum_{j=0}^{\infty}\rho(2^j) < \infty \,.$$

As a matter of fact, by using the subadditivity properties of the conditional expectations of sums Peligrad and Utev [55] showed that the condition

$$\sum_{j=0}^{\infty} \left\| \frac{\mathbb{E}(S_{2^j}|\mathcal{F}_0)}{2^{j/2}} \right\|_2 < \infty \tag{35}$$

is equivalent to (18). Therefore we derived the following corollary:

**Corollary 15** *Assume that a stationary sequence of centered random variables with finite second moments satisfies the condition (31). Then (18) is satisfied and the conclusion of Theorem 11 holds.*

*b) $\alpha$-mixing sequences.*

Strongly mixing sequences are usually treated differently than $\rho$–mixing sequences.

Dedecker and Rio [12] pointed out in their Remark 2 on page 6, that the $\mathbb{L}_2$ condition (16) leads to the suboptimal rates

$$\sum_{k=1}^{\infty} k\int_0^{\tilde{\alpha}(k)} Q^2(u)du \;<\; \infty\,,$$

where $Q$ denotes the cadlag inverse of the function $t \to P(|X_0| > t)$. However, this condition has the advantage that implies that, in the ergodic case the limit of $E(S_n^2)/n$ is strictly positive. For this case Corollary 14 gives the best possible result supported by various counterexamples (see Bradley, [5, Remark 10.23]). By applying Rio's [60] covariance inequality in relation (32), Corollary 14, immediately gives:

**Corollary 16** *Assume*

$$\sum_{k=1}^{\infty} \int_0^{\tilde{\alpha}(k)} Q^2(u)du \; < \; \infty \; . \tag{36}$$

*Then the conclusion of Theorem 13 holds.*

Doukhan, Massart and Rio [17] presented various examples showing that the condition (36) is optimal for the central limit theorem under normalization $\sqrt{n}$. For further discussion and several other optimal features of the above corollary we refer to Bradley [3].

Notice that, as it was pointed out in the above reference, the couple of conditions

$$\mathbb{E}|X_0|^t < \infty, \text{ and } \sum_{k=1}^{\infty} k^{\frac{2}{t-2}} \tilde{\alpha}(k) < \infty, \text{ for some } t > 2 \tag{37}$$

are sufficient for the validity of (36).

Although not being optimal in this case, the general condition (18) when applied to strongly mixing sequences gives results that are surprisingly close to (36). For example, from Corollary 12, the sufficient condition is now

$$\sum_{k=1}^{\infty} \log^v(k) \int_0^{\tilde{\alpha}(k)} Q^2(u)du \; < \; \infty \; ,$$

for some $v > 1$, which gives for the mixing coefficient with the polynomial rate the same sufficient condition as (36).

*c) $\alpha$-mixing sequences plus an extra moment condition.*

In this example, it is shown that the polynomial rates can be improved to logarithmic ones by adding an extra moment condition.

Suppose that $(X_n)_{n\in\mathbb{Z}}$ is a stationary, strongly mixing sequence with $\mathbb{E}(X_0) = 0$ and $\mathbb{E}|X_0|^t < \infty$ for some $t > 2$. Assume that the sequence has a logarithmic mixing rate, $\sum_{r=1}^{\infty} \alpha^{[1/2-1/t]}(2^r) < \infty$ and $\mathrm{Var}(S_n) \to \infty$. Assume also an extra condition:

$$\mathbb{E}|S_n|^t \le K(\mathrm{Var}(S_n))^{t/2}, \text{ for all } n \ge 1 \; . \tag{38}$$

The bound (33) is valid again, (see Bradley and Utev [8, p.118], and Bradley

[5, chapter 8]). So, similarly to (34), write

$$\begin{aligned}
\mathbb{E}[\mathbb{E}(S_n|\mathcal{F}_{-n})^2] &= \operatorname{Cov}(S_n, \mathbb{E}(S_n|\mathcal{F}_{-n})) \\
&\le c\alpha^{1-2/t}(n)(\mathbb{E}|S_n|^t)^{1/t}(\mathbb{E}|\mathbb{E}(S_n|\mathcal{F}_{-n})|^t)^{1/t} \\
&\le c\alpha^{1-2/t}(n)(\mathbb{E}|S_n|^t)^{2/t} \\
&\le c\alpha^{1-2/t}(n)K^{2/t}\|S_n\|_2^2 \\
&\le c_1\alpha^{1-2/t}(n)n\ .
\end{aligned}$$

Now, by using the same arguments following (34), we derive that

$$\|\mathbb{E}(S_{2^{r+1}}|\mathcal{F}_0)\|_2 \le c_1 \sum_{j=0}^{r} 2^{j/2}\alpha^{[1-2/t]/2}(2^j)\,,$$

and by the condition imposed to $\alpha$ we derive

$$\sum_{j=0}^{\infty} \left\| \frac{\mathbb{E}(S_{2^j}|\mathcal{F}_0)}{2^{j/2}} \right\|_2 < \infty\ .$$

This example is motivated by Bradley and Peligrad [7] and Peligrad [49], where sequences with a decomposed strong mixing coefficient are treated. On the other hand, Theorem 13 and Corollary 14 cannot be used to improve the condition (36).

#### *4.2.2. Mixingale-type sequences.*

The mixingales sequences where introduced by McLeish [36, 37]. In the adapted case the mixingale-type conditions are imposed to either $\|\mathbb{E}(X_n|\mathcal{F}_0)\|_1$ or $\|\mathbb{E}(X_n|\mathcal{F}_0)\|_2$, and sometimes to $\|\mathbb{E}(g(X_n)|\mathcal{F}_0)\|_1$, or $\|\mathbb{E}(g(X_n)|\mathcal{F}_0)\|_2$, for $g(x)$ running in certain classes of functions. Notice first that our Corollary 12 and the condition (23), can be viewed as mixingale type results.

Considering the projection operator as defined by (9), we start this section with a general result (see Hannan [23] for the ergodic case, Volný [67] for related results and Corollary 3 in Dedecker and Merlevède [11], for an extension to Hilbert space valued random variables).

**Proposition 17** *Let $(X_i)_{i\in\mathbb{Z}}$ be a stationary sequence with $\mathbb{E}(X_0) = 0$ and $\mathbb{E}(X_0^2) < \infty$ and stationary filtration $(\mathcal{F}_i)_{i\in\mathbb{Z}}$. Let $\mathcal{F}_{-\infty} = \bigcap_{i\in\mathbb{Z}}\mathcal{F}_i$ and define the projection operator $P_i(X) = \mathbb{E}(X|\mathcal{F}_i) - \mathbb{E}(X|\mathcal{F}_{i-1})$. Assume that*

$$\mathbb{E}(X_0|\mathcal{F}_{-\infty}) = 0 \quad a.s. \quad and \quad \sum_{i\ge1}\|P_0(X_i)\|_2 < \infty\ . \tag{39}$$

*Then the conclusion of Theorem 13 holds.*

To comment on the proof of this result, we use the fact that the central limit theorem holds under condition (39) (see for instance Volný [67]). On the other hand, to prove tightness, that is implied by

$$\left(\frac{\max_{1\le i\le n} S_i^2}{n}\right) \quad \text{is uniformly integrable,}$$

the maximal inequality (10) is used. In this maximal inequality the selection of constants that gives the above result is:

$$b_m = 2^{-m-1} \quad \text{and} \quad a_m = (\|P_0(X_m)\|_2 + (m+1)^{-2})^{-1} \; .$$

As a consequence of Proposition (17), the following sufficient condition in the invariance principle for stationary sequences is valid (see condition (1.3) in Dedecker and Merlevède [10] and also remark 6 in Dedecker and Merlevède [11]).

**Corollary 18** *Let $(X_i)_{i\in\mathbb{Z}}$ be a stationary sequence with $\mathbb{E}(X_0)=0$ and $\mathbb{E}(X_0^2)<\infty$. In addition assume that there exists a sequence $(L_k)_{k\ge 1}$ of positive numbers such that*

$$\sum_{k=1}^{\infty}\left(\sum_{i=1}^{k} L_i\right)^{-1} < \infty \quad \textit{and} \quad \sum_{k=1}^{\infty} L_k \|\mathbb{E}(X_k|\mathcal{F}_0)\|_2^2 < \infty \; . \tag{40}$$

*Then (39) is satisfied, and the conclusion of Theorem 13 holds.*

Let us notice that criterion (40) is satisfied if either condition (2.5) in McLeish [39] holds or $(X_i,\mathcal{F}_i)_{i\in\mathbb{Z}}$ is a mixingale of size $-1/2$ (cf. McLeish [37, Definitions (1.2) and (1.4)]). In addition, this criterion is sharp in the sense that the choice $L_k\equiv 1$ is not strong enough to imply the weak convergence of $n^{-1/2}S_n$ (see e.g. Proposition 7 of Dedecker and Merlevède [10]).

By applying an elementary inequality: $(L_1+...+L_n)^{-1}\le n^{-2}(1/L_1+...+1/L_n)$, it is easy to see that

$$\sum_{n\ge 1}(L_1+...+L_n)^{-1} \le C\sum_{n\ge 1}(nL_n)^{-1}$$

and so condition (23) implies (40). Notice however that when $L_n$ is an increasing sequence or $L_n$ is a slowly varying function in the strong sense then (40) and (23) are also equivalent.

To present our next results, we need the following definition.

**Definition 19** *For the integrable random variable $X$, define the "upper tail" quantile function $Q$ by*

$$Q(u) = \inf\{t\ge 0 : \mathbb{P}\left(|X_0|>t\right)\le u\} \; .$$

*Note that, on the set $[0,\mathbb{P}(|X|>0)]$, the function $H: x\to\int_0^x Q(u)du$ is an absolutely continuous and increasing function with values in $[0,\mathbb{E}|X|]$. Denote by $G$ the inverse of $H$.*

**Corollary 20** *Suppose that $(X_n)_{n\geq 0}$ is a stationary sequence with $\mathbb{E}(X_0)=0$ and $\mathbb{E}|X_0|^2<\infty$. Let $Q$ and $G$ be as introduced in Definition 19 for the variable $X_0$. Consider the condition*

$$\sum_{k=1}^{\infty}\int_0^{\|E(X_k|\mathcal{F}_0)\|_1} Q\circ G(u)du \ < \ \infty \,. \tag{41}$$

*We have the implications (36) $\Rightarrow$ (41) $\Rightarrow$ (29) $\Rightarrow$ (28). In particular, under (41), the conclusion of Theorem 13 holds.*

The result follows from Corollary (14), where the series ( 29) is bounded via the following inequality stated in Dedecker and Doukhan [9, Proposition 1]: If $Y$ is $\mathcal{M}$-measurable, then

$$|\,\mathbb{E}(YX)|\leq \mathbb{E}|Y\mathbb{E}(X|\mathcal{M})|\leq \int_0^{\|\mathbb{E}(X|\mathcal{M})\|_1} Q_Y\circ G_X(u)du\,.$$

To make condition (41) more transparent, we present the following result.

**Corollary 21** *Any of the following conditions implies (41).*

1. $\mathbb{P}(|X_0|>x)\leq (c/x)^t$ *for some* $t>2$, *and* $\sum_{k\geq 0}\left(\|\mathbb{E}(X_k|\mathcal{F}_0)\|_1\right)^{(t-2)/(t-1)}<\infty$.
2. $\|X_0\|_t<\infty$ *for some* $t>2$, *and* $\sum_{k\geq 1}k^{1/(t-2)}\|\mathbb{E}(X_k|\mathcal{F}_0)\|_1<\infty$.
3. $\mathbb{E}(X_0^2(\ln(1+|X_0|)))<\infty$ *and* $\|\mathbb{E}(X_k|\mathcal{F}_0)\|_1=O(a^k)$ *for some* $a<1$.

*4.2.3. Linear processes.*

Encountered in various applications such as linear time series, moving averages provide a reach class of examples where different methods of analysis are available. Mixing methods do not appear to provide the best techniques. In order to verify the mixing conditions some additional assumptions are usually imposed to the density of variables.

Applications of martingale technique to linear processes is a traditional technique covered in numerous papers and books specially devoted to linear time series. In this section we briefly analyze causal or one-sided linear processes, i.e. partial sums of the following moving averages,

$$X_k=\sum_{j=-\infty}^{k}a_{k-j}\xi_j=\sum_{i\geq 0}a_i\xi_{k-i}\,,$$

where $\{\xi_k\}_{k\in\mathbb{Z}}$ is a stationary sequence of martingale differences with finite second moment $\sigma^2$, and

$$\sum_{j=0}^{\infty}a_j^2<\infty\,. \tag{42}$$

Since

$$\mathbb{E}(X_k|\mathcal{F}_0) = \sum_{j=-\infty}^{0} a_{k-j}\xi_j = \sum_{i=k}^{\infty} a_i\xi_{k-i} ,$$

we obtain

$$\mathbb{E}|\mathbb{E}(X_k|\mathcal{F}_0)|^2 = \|\mathbb{E}(X_k|\mathcal{F}_0)\|_2^2 = \sigma^2 \sum_{i=k}^{\infty} a_i^2 \to 0 \text{ as } k \to \infty .$$

Therefore the mixing-type condition (30) is satisfied.

The moving average is usually called long range dependent if

$$\sum_{j=0}^{\infty} |a_j| = \infty , \tag{43}$$

and the linear process is called short range dependent if

$$\sum_{j=0}^{\infty} |a_j| < \infty . \tag{44}$$

The central limit theorem and its invariance principle for short range dependent variables are known for a while (see, e.g. Hannan [23]). The central limit theorem under condition (44) follows easily from the following representation

$$\sum_{k=1}^{n} X_k - A\sum_{j=1}^{n} \xi_j = \sum_{j=-\infty}^{\infty} \left(\sum_{k=1}^{n} c_{k-j}\right) \xi_j$$

where

$$A = \sum_{j=0}^{\infty} a_j \ , \ c_0 = a_0 - A \quad \text{and} \quad c_i = a_i \quad \text{for} \quad i \geq 1 \ \text{ and } \ c_i = 0 \quad \text{for} \quad i \leq -1 \ .$$

This representation implies

$$\frac{1}{n}\mathbb{E}\left|\sum_{k=1}^{n} X_k - A\sum_{j=1}^{n} \xi_j\right|^2 = \mathbb{E}[\xi^2]\frac{1}{n}\sum_{j=-\infty}^{\infty}\left|\sum_{i=1-j}^{n-j} c_i\right|^2 \to 0 \quad (n \to \infty)$$

Because $\sum_{i=0}^{\infty} c_i = 0$, by condition (44), and Lemma 1 in Merlevède, Peligrad and Utev [43] the last term in the above inequality is convergent to 0 as $n \to \infty$ and the CLT follows for $S_n = \sum_{k=1}^{n} X_k$ with the normalization $\sqrt{n}$ by the standard CLT for martingales. A similar approach was used in Merlevède, Peligrad and Utev [43] for example, applied for the case when the innovations are dependent or/and Hilbert space valued random variables and $(a_j)_{j\geq 0}$ is a sequence of operators.

In a recent paper, Peligrad and Utev [56] showed that, as a matter of fact, despite the long range dependence, the central limit theorem holds for partial sums of a linear process when the sequence of constants satisfies condition (42) and the sequence of stationary innovations $\{\xi_k\}_{k\in\mathbb{Z}}$ satisfies(39) (for the martingale difference innovations, the result is formulated in Proposition 31 below, in part 4.4.1.)

The invariance principle is more involved and in general does not hold under just the condition (42) (see for instance, the example 1, Section 3.2 in Merlevède and Peligrad [42]).

Now, we shall analyze the implications of various martingale results described in the previous sections on the invariance principle for linear processes.

Under condition (44), both the central limit theorem and invariance principle follow from the projective criteria that is Proposition 17. Moreover, for the causal linear processes, condition (39) is actually equivalent to (44), which is easily seen from the following well-known argument (Hall and Heyde, [22])

$$P_0(X_k) = \sum_{j=-\infty}^{k} a_{k-j} P_0(\xi_j) \ = \ \sum_{j=-\infty}^{0} a_{k-j}\xi_j - \sum_{j=-\infty}^{-1} a_{k-j}\xi_j = a_k \xi_0 \ .$$

Thus, $\|P_0(X_k)\|_2 = \|\xi_0\|_2 |a_k|$ and so in this case

$$\sum_{i\in\mathbb{Z}} \|P_0(X_i)\|_2 < \infty \quad \text{if and only if} \quad \sum_{i=0}^{\infty} |a_i| < \infty$$

It is easy to see that

$$\mathbb{E}(S_n|\mathcal{F}_0) = \sum_{j=-\infty}^{0} \xi_j \left( \sum_{k=1}^{n} a_{k-j} \right) \ .$$

We then immediately notice that if $\xi_1$ has a normal distribution, then conditions (16) and (17) are equivalent.

Wu (2002, Lemma 1), pointed out that the mixingale type condition from Corollary 12, applied to the one-sided moving averages is equivalent to

$$\sum_{n=1}^{\infty} \frac{1}{\sqrt{n}} \left( \sum_{k=n}^{\infty} a_k^2 \right)^{1/2} < \infty$$

and it is slightly stronger than (44). Moreover, Wu proved (2002, Proposition 1), that

$$\sup_n \|\mathbb{E}(S_n|\mathcal{F}_0)\|_2 < \infty \tag{45}$$

if and only if $\sum_{i\geq 0} a_i$ exists and

$$\sum_{n=1}^{\infty} \left( \sum_{k=n}^{\infty} a_k \right)^2 < \infty \ . \tag{46}$$

The relation in above first appears in Hall and Heyde [?, p.146]. In fact, (46) is the necessary and sufficient for the following coboundary representation

$$X_k = A\xi_k + Y_{k-1} - Y_k$$

where

$$A = \sum_{i=0}^{\infty} a_i \ , \ Y_k = \sum_{j=0}^{\infty} \hat{a}_j \xi_{k-j} \ , \ \hat{a}_j = \sum_{i=j+1}^{\infty} a_i .$$

A sufficient conditions for (46) is $\Sigma_{j=1}^{\infty} j^{1/2}|a_j| < \infty$ as it was pointed out in Phillips and Solo [59], who proved several limit theorems, including the central limit theorem, invariance principle and the law of the iterated logarithm under this and similar conditions.

Examples of long range dependent moving averages which also satisfy condition (46) were known since Heyde [26], who suggested to use the spectral density approach in the analysis of the Gordin type condition (16). Some related examples are suggested in Propositions 3 and 4 in Wu [69]. Thus, even the stronger condition (45) than condition (18) is satisfied for some long range dependent sequences.

The following result is an immediate corollary of Theorem 11.

**Corollary 22** *Let $b_n = a_0 + \cdots + a_n$ and assume that*

$$\sum_{n=1}^{\infty} \frac{1}{n^{3/2}} \Big\{ \sum_{j=0}^{\infty} (b_{n+j} - b_j)^2 \Big\}^{1/2} \ < \ \infty \ . \tag{47}$$

*Then (18) is satisfied and the conclusion of Theorem 11 holds.*

As we have seen the conditions (44), (46), and (47) all assure the validity of the invariance principle. Notice that (47) is in general weaker than condition (46) (just compare their equivalent (45) and (18); in addition Example 2 in part 4.4.1 below shows that (47) does not imply (46)). However conditions (44) and (47) are independent, one can hold without the other one to hold.

#### *4.2.4. Shift processes.*

Maxwell and Woodroofe [35] have specialized their central limit theorem and invariance principle to one-sided shift processes, also known under the name of Raikov sums. We will consider only Bernoulli shifts and show what kind of improvement can be made by using Theorem 11. Let $\{\varepsilon_k; k \in \mathbb{Z}\}$ be an i.i.d. sequence of random variables with $\mathbb{P}(\varepsilon_1 = 0) = \mathbb{P}(\varepsilon_1 = 1) = 1/2$ and let

$$Y_n = \sum_{k=0}^{\infty} 2^{-k-1} \varepsilon_{n-k} \quad \text{and} \quad X_n = g(Y_n) - \int_0^1 g(x)dx \, ,$$

where $g \in \mathbb{L}_2(0,1)$, $(0,1)$ being equipped with the Lebesgue measure.

The following result is an immediate consequence of Proposition 3 in Maxwell and Woodroofe [35] and of Theorem 11.

**Corollary 23** *For the Bernoulli shift process, if* $g \in \mathbb{L}_2(0,1)$*, and*

$$\int_0^1 \int_0^1 [g(x) - g(y)]^2 \frac{1}{|x-y|} (\log[\log \frac{1}{|x-y|}])^t dxdy \ < \ \infty \tag{48}$$

*for some* $t > 1$*, then (24) is satisfied and the conclusion of Theorem 11 holds.*

As a concrete example of a map which is not covered by Maxwell and Woodroofe's theorem, we consider the function

$$g(x) = \frac{1}{\sqrt{x}} \frac{1}{[1 + \log(2/x)]^4} \sin(\frac{1}{x}) \ , \ 0 < x < 1 \ .$$

The analysis of this example is based on the idea of proof developed by Maxwell and Woodroofe [35] who treated the function

$$g_\alpha(x) = \frac{1}{x^\alpha} \sin(\frac{1}{x}) \ , \ 0 < x < 1,$$

where $0 < \alpha < 1/2$ (actually, $g_\alpha \in \mathbb{L}_p(0,1)$ for all $p < 1/\alpha$).

To prove the convergence of the integral (48), we shall consider the case $t = 2$ and change the variables from $x$ to $1/x$ and from $y$ to $1/y$, and introduce the function

$$G(x) = \sqrt{x} \sin(x) [1 + \log 2x]^{-4} \ .$$

In addition, we split the resulting integral into the following three regions

$$D_1 = \{x > 1, y > 1, |x-y| > 1\} \ , D_2 = \{x > 1, y > 1, |x-y| \le 1; x + y \ge 5\};$$
$$D_3 = \{x > 1, y > 1, |x-y| \le 1, x + y \le 5\} \ .$$

and so we have to show that for $i = 1, 2, 3$,

$$J_i = \int\int_{D_i} [G(x) - G(y)]^2 \frac{1}{xy|x-y|} (\log|\log \frac{xy}{|x-y|}|)^2 dxdy \ < \ \infty \ .$$

We first notice that $J_3 < \infty$ because the region of integration is finite and the integrand is uniformly bounded since the function $G$ is continuously differentiable for bounded $x$.

We are going to show that $J_1 < \infty$. The analysis of $J_2$ is similar.

By using $(a + b) \le 2a^2 + 2b^2$ and symmetry, $J_1$ is bounded up to a factor 4 by

$$\begin{aligned}
&\int\int_{x>1, y>1, |x-y|>1} [G(x)]^2 \frac{1}{xy|x-y|} (\log|\log \frac{xy}{|x-y|}|)^2 dxdy \\
&\quad \le \int\int_{x>3, y>1, y \le 2x, |y-x|>1} \frac{1}{[1 + \log 2x]^8} \frac{1}{y|x-y|} (\log|\log \frac{xy}{|x-y|}|)^2 dxdy \\
&\quad + \int\int_{x>1, y>1, y \ge 2x} \frac{1}{[1 + \log 2x]^8} \frac{1}{y|x-y|} (\log|\log \frac{xy}{|x-y|}|)^2 dxdy \\
&\quad + \int\int_{1<x<3, y>1, y \le 2x, |y-x|>1} \frac{1}{[1 + \log 2x]^8} \frac{1}{y|x-y|} (\log|\log \frac{xy}{|x-y|}|)^2 dxdy \\
&= I_1 + I_2 + I_3 \ .
\end{aligned}$$

Again, it is easy to see that $I_3 < \infty$, since the region of integration is finite and the integrand is uniformly bounded. And, we only show that $I_1 < \infty$ with $I_2$ being treated in the similar way.

We notice that for $x > 3, y > 1, y \leq 2x, |y-x| > 1$, we have $x/(x-1) \leq xy/|x-y| \leq 2x^2$, and so $1/2x \leq \log \frac{xy}{|x-y|} \leq 2x$ which implies

$$\left|\log\log \frac{xy}{|x-y|}\right| \leq \log 2x \ .$$

Consequently, by the Fubini formula,

$$\begin{aligned} I_1 \quad \leq \quad & \int_{x=3}^{\infty} \frac{1}{[1+\log 2x]^6}\Big(\int_{y=x+1}^{y=2x} \frac{1}{y(y-x)} dy\Big) dx \\ & + \int_{x=3}^{\infty} \frac{1}{[1+\log 2x]^6}\Big(\int_{y=1}^{y=x-1} \frac{1}{y(x-y)} dy\Big) dx < \infty \ . \end{aligned}$$

*4.2.5. A particular Markov chain example*

The following example is motivated by Isola [32], who applied discretization to the Pomeau–Manneville type 1 intermittency model. Let $\{Y_k; k \geq 0\}$ be a discrete Markov chain with the state space $Z^+$ and transition matrix $P = (p_{ij})$ given by $p_{k(k-1)} = 1$ for $k \geq 1$ and $p_j = p_{0(j-1)} = P(\tau = j)$, $j = 1, 2, \ldots$, (that is whenever the chain hits 0, $Y_t = 0$, it then regenerates with the probability $p_j$). Now, unlike Peligrad and Utev [55], we do not ask that $p_1, p_2 > 0$. We only need that $p_{n_j} > 0$ along $n_j \to \infty$. We observe that if from state 0 we can go to state $j$ then, we can go from state 0 to states $j-1, \ldots, 0$. That is the chain is then irreducible. The stationary distribution exists if and only if $\mathbb{E}[\tau] < \infty$ and it is given by

$$\pi_j = \pi_0 \sum_{i=j+1}^{\infty} p_i \ , \ j = 1, 2 \ldots$$

where $\pi_0 = 1/\mathbb{E}[\tau]$.

Further, as in Peligrad and Utev [55], as a functional $g$ we take $\mathbb{I}_{(x=0)} - \pi_0$ , where $\pi_0 = \mathbb{P}_\pi(Y_0 = 0)$ under the stationary distribution denoted by $\mathbb{P}_\pi$ ( $\mathbb{E}_\pi$ denotes the expectations for the process started with the stationary distribution). The stationary sequence is then defined by

$$X_j = \mathbb{I}_{(Y_j=0)} - \pi_0 \quad \text{so that} \quad S_n = \sum_{j=1}^{n} X_j = \sum_{j=1}^{n} \mathbb{I}_{(Y_j=0)} - n\pi_0 \ .$$

Finally, a particular distribution $(p_k)$ was constructed in Peligrad and Utev [55] to show the optimality of the Theorem 11 in the following form.

**Theorem 24** *For any non-negative sequence $a_n \to 0$ there exists a stationary ergodic discrete Markov chain $(Y_k)_{k\geq 0}$ and a functional $g$ such that $X_i = g(Y_i)$;*

$i \geq 0$, $\mathbb{E}[X_1] = 0$, $\mathbb{E}[X_1^2] < \infty$ *and*

$$\sum_{n=1}^{\infty} a_n \frac{\|\mathbb{E}(S_n|Y_0)\|_2}{n^{3/2}} < \infty \text{ but } \frac{S_n}{\sqrt{n}} \text{ is not stochastically bounded}\,. \tag{49}$$

It is interesting to notice that if $p_{2k} = 0$ for all integer $k$, then the chain is periodic with period 2 and so, because of the lack of mixing property, the condition (30) is not valid in this situation. On the other hand, a similar analysis to the one in Peligrad and Utev [55] shows that condition (18) is satisfied if $\mathbb{E}[\tau^2 \log(1+\tau)^4] < \infty$. This condition is close to the optimal condition for CLT that follows from the central limit theorem for renewal sequences, namely $\mathbb{E}[\tau^2] < \infty$.

### 4.3. Towards a more general normalization

In Section 4.1, we have surveyed several aspects of the central limit theorem and invariance principle under conditions assuring that the variance of partial sums is linear in $n$. More precisely, we treated the case when

$$\lim_{n\to\infty} \frac{\operatorname{Var}(S_n)}{n} = \sigma^2 < \infty\,.$$

In this part, we consider the central limit theorem and invariance principle under more general normalization such as

$$\sigma_n = \sqrt{\operatorname{Var}(S_n)} = \operatorname{Stdev}(S_n) \quad \text{or} \quad b_n = \sqrt{\frac{\pi}{2}}\ \mathbb{E}|S_n|\,.$$

The first normalization was used in many papers involving mixing structures probably initiated by Ibragimov [27, 30]. The central limit theorem and invariance principle under the second normalization were considered in Dehling, Denker and Philipp [13] and Peligrad [48], respectively.

When the normalization is $\sigma_n = \operatorname{Stdev}(S_n)$ and $S_{[nt]}/\sigma_n$ converges in distribution to the Brownian motion, then, necessarily $\sigma_n^2$ has the representation $\sigma_n^2 = nh(n)$ with $h(n)$ a slowly varying function at infinity. It is natural then to give sufficient conditions for this representation of $\sigma_n^2$.

**Proposition 25** *Assume that $(X_i)_{i\in\mathbb{Z}}$ is a stationary sequence with $\mathbb{E}(X_0) = 0$ and $\mathbb{E}(X_0^2) < \infty$. Assume that $\sigma_n^2 \to \infty$, as $n \to \infty$, and that one of the following conditions is satisfied:*

*(i)*

$$\|\mathbb{E}(S_n|\mathcal{F}_0)\|_2 = o(\sigma_n) \text{ as } n \to \infty\,; \tag{50}$$

*(ii)* $\lim_{n\to\infty} q_0(n) = 0$, *where*

$$q_0(n) = \sup\{|\operatorname{Cov}(S_k, S_{a+k} - S_a)|/\|S_k\|_2^2\ ;\ k \geq 1, a \geq n+k\}$$

*(iii)*

$$\liminf_{n\to\infty} \frac{\operatorname{Var}(S_n)}{n} > 0 \quad \text{and} \quad n\mathbb{E}(X_0 X_n) \to 0 \text{ as } n \to \infty\,.$$

*Then $\sigma_n^2 = nh(n)$ with $h(n)$ a slowly varying function in the strong sense at infinity.*

The statement (i) was proved by Wu and Woodroofe [70]. The result (ii) is due to Ibragimov [27, 28, 29]; the proof can be found in Bradley [5, Theorem 8.13]. The last statement was proved in Merlevède and Peligrad [41].

The first result we would like to mention is in the spirit of Theorem 1 in Dedecker and Merlevède [10] and it is a version of Theorem 1 in Merlevède and Peligrad [42].

**Proposition 26** *Let $(X_i)_{i\in\mathbb{Z}}$ be a stationary sequence with $\mathbb{E}(X_0) = 0$ and $\mathbb{E}(X_0^2) < \infty$. Set $S_n := \sum_{k=1}^n X_k$ and $\sigma_n^2 = \mathrm{Var}(S_n)$. Assume that $\sigma_n^2 \to \infty$, as $n \to \infty$, and that*

$$\frac{S_n^2}{\sigma_n^2} \text{ is an uniformly integrable family} \tag{51}$$

$$\|\mathbb{E}(S_n|\mathcal{F}_0)\|_1 = o(\sigma_n) \text{ as } n \to \infty, \tag{52}$$

*and that there exists a positive random variable $\eta$ such that*

$$\lim_{n\to\infty} \mathbb{E}(\frac{S_n^2}{\sigma_n^2}|\mathcal{F}_{-n}) = \eta \text{ in } \mathbb{L}_1 \ . \tag{53}$$

*Then the random variable $\eta$ is $\mathcal{I}$ – measurable and $\sigma_n^{-1}S_n \stackrel{\mathcal{D}}{\to} \sqrt{\eta}\,N$, where $N$ is a standard Gaussian r.v. independent of $\mathcal{I}$.*

Notice that (51) together with (52) imply (50), and then, according to Proposition 25, the fact that $\sigma_n^2 = nh(n)$ with $h(n)$ a slowly varying function in the strong sense at infinity.

An additional condition is required for the invariance principle:

**Proposition 27** *Under conditions of Proposition 26, assume in addition that*

$$\lim_{\lambda\to\infty} \limsup_{n\to\infty} \lambda^2 \mathbb{P}\left(\max_{1\le i\le n} \frac{|S_i|}{\sigma_n} \ge \lambda\right) = 0 \ . \tag{54}$$

*Then $\sigma_n^{-1}W_n \stackrel{\mathcal{D}}{\Longrightarrow} \sqrt{\eta}W$ , where $W$ is a standard Brownian motion independent of $\mathcal{I}$ and $\eta$ is the $\mathcal{I}$-measurable non-negative variable, defined in Proposition 26.*

Clearly conditions (51) and (54) are satisfied if

$$\frac{\max_{1\le i\le n} S_i^2}{\sigma_n^2} \text{ is uniformly integrable} \ . \tag{55}$$

The following corollary will be applied to analyze examples mentioned in part 4.4.

**Corollary 28** *Assume that $\sigma_n^2 \to \infty$, as $n \to \infty$, and that the conditions (38), (52) and (53) are satisfied. Then, the conclusion of Proposition 27 holds.*

**Proof**. Notice first that clearly (38) implies (51) which, as we mentioned before, combined with (52) implies that $\mathrm{Var}(S_n) = \sigma_n^2 = nh(n)$, with $h(n)$ a slowly varying function in the strong sense at infinity. Then, by applying Lemma 3.7 in Peligrad [46] with $t$ instead of 4, it follows from (38) and the fact that $\sigma_n^2 = nh(n)$, with $h(n)$ a slowly varying function, that

$$\mathbb{E}\Big(\max_{1\le i\le n} |S_i|^t\Big) \le c_t \mathrm{Var}(S_n)^{t/2} \, .$$

which implies (55) and the corollary now follows from Proposition 27.

By using Proposition 27 it is possible to extend Theorem 13 of Dedecker and Rio [12]. The following result is in this direction (see Corollary 1 in Merlevède and Peligrad [42]):

**Corollary 29** *Let $(X_i)_{i\in\mathbb{Z}}$ be a stationary sequence with $\mathbb{E}(X_0) = 0$ and $\mathbb{E}(X_0^2) < \infty$. In addition assume that $\sigma_n^2 = nh(n)$ where $h(n)$ is slowly varying in the strong sense and satisfies $\lim_{n\to\infty} h(n) = c$, where $c$ is a strictly positive constant possibly infinite, and that*

$$\frac{n}{\sigma_n^2}\mathbb{E}(X_0 S_n|\mathcal{F}_0) \text{ is convergent in } \mathbb{L}^1 \text{ to a random variable } \mu \text{ as } n \to \infty \, . \quad (56)$$

*Then the invariance principle of Proposition 27 holds with*

$$\eta = c^{-1}\mathbb{E}(X_0^2|\mathcal{I}) + 2\mathbb{E}(\mu|\mathcal{I}) \, .$$

The proofs were based on the Bernstein blocking techniques of partitioning the variables in big and small blocks of suitable sizes combined with the martingale decomposition of the sums of variables in blocks. More precisely the variables are divided in blocks of size $p = o(n)$ followed by smaller blocks of size $q = o(p)$. Then the sum of variables in the big blocks (of size $p$), denoted by $Y_j, 1 \le j \le k = [n/(p+q)]$ is decomposed into a martingale

$$M_k = \sum_{j=1}^{k}(Y_j - \mathbb{E}(Y_j|Y_1 ... Y_{j-1}))$$

and the error term

$$R_k = \sum_{j=1}^{k} \mathbb{E}(Y_j|Y_1 ... Y_{j-1}) \, .$$

A careful analysis of the non stationary martingale and the error term leads to sufficient conditions for the convergence of the sequence $b_n^{-1}W_n \overset{\mathcal{D}}{\Longrightarrow} W$, where $W$ is a standard Brownian motion and $b_n = \sqrt{\frac{\pi}{2}}\ \mathbb{E}|S_n|$ .

In a recent paper Merlevède and Peligrad [42] obtained a general invariance principle under the above normalization involving the absolute first moment of partial sums. This normalization leads to more general sufficient conditions for the invariance principle. We shall give here, for simplicity, only a corollary that was derived from Theorem 3 and Corollary 2 in Merlevède and Peligrad [42] under the additional assumption:

$$\liminf_{n\to\infty} \frac{\mathrm{Var}(S_n)}{n} > 0 . \tag{57}$$

**Proposition 30** *Assume that $\sigma_n^2 = nh(n)$ with $h(n)$ a slowly varying function in the strong sense as $n \to \infty$ and that the conditions (53) and (57) are satisfied. Moreover suppose that*

**(i)** $\mathbb{P}(|X_0| \le T) = 1$ *for a positive number $T$ and* $\lim_{n\to\infty} \|\mathbb{E}(S_n|\mathcal{M}_{-n})\|_1 = 0$,

*or that*

**(ii)** *there exist $r > 2$ and $c > 0$ such that $\mathbb{P}(|X_0| > x) \le (c/x)^r$ and* $\lim_{n\to\infty} n^{1/(r-1)} \big(\|\mathbb{E}(S_n|\mathcal{M}_{-n})\|_1\big)^{(r-2)/(r-1)} = 0$.

*Then*

$$\frac{W_n(t)}{\sqrt{\frac{\pi}{2}}\mathbb{E}|S_n|} \overset{\mathcal{D}}{\Longrightarrow} \sqrt{\eta} W \ \text{ as } n \to \infty \ ,$$

*where $W$ is a standard Brownian motion independent of $\mathcal{I}$ and $\eta$ is the $\mathcal{I}$-measurable non-negative variable, defined in Proposition 26.*

Merlevède and Peligrad [42] obtained sharp sufficient conditions for the invariance principle when

$$\lim\inf_{n\to\infty} \frac{\mathrm{Var}(S_n)}{n} = 0. \tag{58}$$

The results of this type are more delicate and further study is needed, since there are several examples of dependent classes of random variables satisfying (58) such as suggested in Ibragimov and Rozanov [31], Herrndorf [24], Bradley [5, v.3, 27.5 and 27.6], Wu and Woodroofe [70] and in Merlevède and Peligrad [42].

### *4.4. Examples. Part 2*

#### *4.4.1. Linear processes and general normalization.*

In this section we present an invariance principle for a class of linear processes when the variance of partial sums is not necessarily linear in $n$. This part was motivated by Wu and Woodroofe [70] and Merlevède and Peligrad [42].

Let $(\xi_i)_{i\in\mathbb{Z}}$ be a stationary sequence of martingale differences with finite second moment $\sigma^2$. In addition, as before, let $(a_i, i \ge 0)$ be a sequence of real numbers such that $\sum_{i\ge 0} a_i^2 < \infty$.

As in part 4.2.3, we consider the causal linear process defined by

$$X_k = \sum_{i\geq 0} a_i \xi_{k-i} \quad \text{and let} \quad b_n = a_0 + \cdots + a_n \ . \tag{59}$$

Notice that

$$\sigma_n^2 = \mathrm{Var}(S_n) = \sigma^2 \sum_{k=0}^{n-1} b_k^2 \ + \ \sigma^2 \sum_{j=0}^{\infty} (b_{n+j} - b_j)^2 \ .$$

By using the representation of $S_n$ as a non-stationary martingale

$$S_n = \sum_{j\in\mathbb{Z}} \xi_j \Big( \sum_{k=1}^{n} a_{k-j} \Big)$$

and applying the martingale approximation method combined with the Bernstein blocking technique, Peligrad and Utev [56] proved the following central limit theorem, and further reaching results.

**Proposition 31** *Assume $(X_k)_{k\in\mathbb{Z}}$ is a linear process defined as above and suppose that $\sigma_n^2 = \mathrm{Var}(S_n) \to \infty$. Then, $\sigma_n^{-1} S_n \xrightarrow{\mathcal{D}} \sqrt{\sigma^{-2}\mathbb{E}(\xi_0^2|\mathcal{I})}\, N$, where $N$ is a standard Gaussian r.v. independent of $\mathcal{I}$.*

Let us now treat $S_n$ as a partial sum of a stationary sequence. As it follows from part 4.2.3,

$$\|\mathbb{E}(S_n|\mathcal{F}_0)\|_2^2 \ = \ \sigma^2 \sum_{j=0}^{\infty} (b_{n+j} - b_j)^2 \ .$$

Using this fact, Wu and Woodroofe [70] noticed that the necessary and sufficient condition for (50), i.e $\|\mathbb{E}(S_n|\mathcal{F}_0)\|_2 = o(\sigma_n)$ together with $\sigma_n \to \infty$ is

$$\sum_{k=0}^{n-1} b_k^2 \to \infty, \quad \text{as} \ \ n \to \infty, \tag{60}$$

and

$$\sum_{j=0}^{\infty} (b_{n+j} - b_j)^2 = o\big(\sum_{k=0}^{n-1} b_k^2\big) \ . \tag{61}$$

Hence, under conditions (60) and (61), the central limit theorem follows from Proposition 26. This traditional approach of proving the central limit theorem for linear processes via general results for stationary sequences did not appear to be as efficient as Proposition 31. However, it has two advantages. First, it allows to find checkable sufficient conditions for the invariance principle. Second, it allows to establish a so-called conditional central limit theorem, that is roughly speaking a stable convergence with respect to the starting $\sigma$–field $\mathcal{F}_0$. In both cases, $\sigma_n^2 = nh(n)$ is a necessary condition. Moreover, Wu and Woodroofe [70,

Lemma 2] proved that the stronger condition (50) is necessary for the conditional central limit theorem.

We shall mention the following two examples of linear processes that were analyzed by Wu and Woodroofe [70] and that satisfy the conditions ( 60) and (61).

*Example 1.* We consider first $\mathbb{E}(\xi_0)^2 = 1$ and the selection of constants

$$a_0 = 0 \quad \text{and} \quad a_n = \frac{1}{n} \quad \text{for} \quad n \geq 1 \,.$$

Then, we derive

$$b_n \sim \log(n), \quad \sum_{k=0}^{n-1} b_k^2 \sim n \log^2(n), \text{ and } \sum_{j=0}^{\infty} (b_{n+j} - b_j)^2 = O(n)$$

so the condition (61) holds.

*Example 2.* We consider $\mathbb{E}(\xi_0)^2 = 1$ and the selection of constants

$$a_0 = 0 \quad \text{and} \quad a_1 = \frac{1}{\log(2)}, \quad \text{and} \quad a_n = \frac{1}{\log(n+1)} - \frac{1}{\log(n)} \quad \text{for} \quad n \geq 2 \,.$$

Then, we derive

$$\sum_{k=0}^{n-1} b_k^2 \sim \frac{n}{\log^2(n)}, \text{ and } \sum_{j=0}^{\infty} (b_{n+j} - b_j)^2 = O\left(\frac{n}{\log^3(n)}\right)$$

so, the condition (61) holds. Notice that for this example : $\lim_{n\to\infty} \mathbb{E}(S_n)^2/n = 0$. And although condition (47) holds, Corollary 22 only implies the convergence to the degenerate law.

Based on the second example considered by Wu and Woodroofe [70] that was further developed in Merlevède and Peligrad [42] we can see that even conditions (60) and (61) are not sufficient to ensure that the linear processes $(X_k)_{k\in\mathbb{Z}}$ satisfies the invariance principle:

**Proposition 32** *There is a linear process $(X_k)_{k\in\mathbb{Z}}$ satisfying the central limit theorem and conditions (60) and (61), but such that the invariance principle does not hold.*

In order to derive the weak invariance principle for the class of linear processes studied in this section, we impose an additional condition on the moment of the random variables $(\xi_i)_{i\in\mathbb{Z}}$. The proof of the next theorem is based on Corollary 28 (see Merlevède and Peligrad [42]).

**Proposition 33** *Let $(\xi_i)_{i\in\mathbb{Z}}$ be a stationary sequence of martingale differences such that $\mathbb{E}|\xi_0|^t < \infty$ for some $t > 2$. Let $(a_i)_{i\geq 0}$ be a sequence of real numbers such that $\sum_{i\geq 0} a_i^2 < \infty$ and satisfying (60) and (61). Let $(X_i)_{i\in\mathbb{Z}}$ be the linear process defined by (59). Then, $\sigma_n^{-1} W_n \stackrel{\mathcal{D}}{\Longrightarrow} \sqrt{\sigma^{-2}\mathbb{E}(\xi_0^2|\mathcal{I})}\, W$ where $W$ is a standard Brownian motion independent of $\mathcal{I}$.*

*4.4.2. Strongly mixing sequences and general normalization*

We would like to start this section by a proposition related to the condition (38). For strongly mixing sequences, provided that $\sigma_n^2 = \mathrm{Var}(S_n) \to \infty$, condition (38) alone is sufficient for the invariance principle under normalization $\sigma_n$.

**Proposition 34** *Assume that $(X_i)_{i\in\mathbb{Z}}$ is a stationary, strongly mixing sequence with $\mathbb{E}(X_0) = 0$ and $\mathbb{E}(|X_0|^t) < \infty$ , for some $t > 2$, and $\sigma_n^2 \to \infty$. Assume that (38) is satisfied. Then, $\sigma_n^{-1} W_n \overset{\mathcal{D}}{\Longrightarrow} W$.*

**Proof**. We show that conditions of Corollary 28 are satisfied. We first notice that, by stationarity, for all $m = 1, 2, \ldots$, we have:

$$\mathbb{E}|\,\mathbb{E}(S_n|\mathcal{F}_0)| \le \mathbb{E}|\,\mathbb{E}(S_{n+m}|\mathcal{F}_0)| + \sigma_m \le \mathbb{E}|\,\mathbb{E}(S_n|\mathcal{F}_{-m})| + 2\sigma_m \ .$$

Since $\sigma_n^2 \to \infty$, (52) will follow if we can prove that

$$\lim_{m\to\infty} \lim_{n\to\infty} \frac{\mathbb{E}|\,\mathbb{E}(S_n|\mathcal{F}_{-m})|}{\sigma_n} = 0 \ . \tag{62}$$

But, similarly to the point c) in part 4.2.1., from (38), we derive that there exist positive constants $c_1, c_2$ such that for all $m = 1, 2, \ldots$,

$$\mathbb{E}(\mathbb{E}(S_n|\mathcal{F}_{-m}))^2 \le c_1 \alpha(m)^{1-2/t} (\mathbb{E}|S_n|^t)^{2/t} \le c_2 \alpha(m)^{1-2/t} \sigma_n^2$$

which implies (62) by using the facts that $\sigma_n^2 \to \infty$ and $\alpha(m) \to 0$.

To prove (53), we apply (32) with $X = 1$, $Y = S_n^2/\sigma_n^2 - 1$ and $\mathcal{F} = \mathcal{F}_{-n}$. Then, we successively derive that, for some positive constants $c_3, c_4, c_5$,

$$\begin{aligned}
&\mathbb{E}|\,\mathbb{E}(S_n^2/\sigma_n^2 - 1|\mathcal{F}_{-n})| = \mathrm{Cov}((\ \mathbb{1}_{\mathbb{E}(Y|\mathcal{F}_{-n})>0} - \ \mathbb{1}_{\mathbb{E}(Y|\mathcal{F}_{-n})\le 0}), Y) \\
&\le c_3 \alpha^{1-2/t}(n) [\mathbb{E}(|Y|^{t/2})]^{2/t} = c_3 \alpha^{1-2/t}(n) [\mathbb{E}(|S_n^2/\sigma_n^2 - 1|^{t/2})]^{2/t} \\
&\le c_4 \alpha(n)^{1-2/t} (\mathbb{E}|S_n|^t)^{2/t} \le c_5 \alpha(n)^{1-2/t} \sigma_n^2
\end{aligned}$$

which implies (53) with $\eta = 1$. The proof is now complete.

We shall finish this section by stating the invariance principle for strongly mixing sequences under the normalization

$$b_n = \sqrt{\frac{\pi}{2}}\ \mathbb{E}|S_n| \ .$$

which was established in Merlevède and Peligrad [41].

**Theorem 35** *Suppose that $\{X_k, k \in \mathbb{Z}\}$ is a stationary, strongly mixing sequence with $\mathbb{E}(X_0) = 0$ and $\mathbb{E}(X_0^2) < \infty$. Assume that*

$$\liminf_{n\to\infty} \frac{\mathbb{E}(S_n^2)}{n} > 0 \ . \tag{63}$$

*and*

$$\int_0^{\alpha_n} Q^2(u)du \;=\; o\Big(\frac{1}{n}\Big) \qquad \text{as } n \to \infty\ .$$

*Then,* $b_n^{-1} W_n \stackrel{\mathcal{D}}{\Longrightarrow} W$.

Notice that the normalization cannot be just $\sqrt{n}$. Moreover, the result is optimal in the view of the counterexample by Bradley [3]. This theorem was extended by Merlevède [40] to Hilbert space valued random variables.

The following corollary extends the corresponding results of Ibragimov [27].

**Corollary 36** *Suppose that* $(X_k)_{k\in\mathbb{Z}}$ *is a stationary, strongly mixing sequence with* $\mathbb{E}(X_0) = 0$, $\mathbb{E}(X_0^2) < \infty$ *and satisfying (63). In addition assume that one of the following conditions is satisfied*
*1)* $\mathbb{E}(|X_0|^t) < \infty$ *for some* $t > 2$, *and* $n\alpha_n^{\frac{t-2}{t}} \to 0$ *as* $n \to \infty$.
*2)* $X$ *is bounded and* $n\alpha_n \to 0$ *as* $n \to \infty$.
*Then the conclusion of Theorem 35 holds.*

**Acknowledgment.** The authors thank the anonymous referee for helpful comments and suggestions.